\documentclass[11pt,a4paper]{article}
\usepackage{amsmath,amsfonts,amssymb}
\usepackage[english]{babel}

\newenvironment{proof}{\noindent {\em Proof.}}{ \hfill $\Box$ \\}

\newtheorem{lem}{Lemma}
\newtheorem{pro}[lem]{Proposition}
\newtheorem{rem}[lem]{Remmark}
\newtheorem{ex}[lem]{Example}

\begin{document}

\begin{center}
{\Large{A note on factorisation of division polynomials}}\\
\vspace*{.5cm}
\end{center}
\begin{center} D. Sadornil\\
\vspace*{.5cm} Departamento of Matem\'aticas,U. Salamanca\\ Pza. de
la Merced 1-4,
37008 Salamanca, Spain.\\ {\small email: {\tt sadornil@usal.es}}\\
\end{center}

\begin{abstract}
In \cite{Verdure}, Verdure gives the factorisation patterns of
division polynomials of elliptic curves defined over a finite field.
However, the result given there contains a mistake. In this paper,
we correct it.
\end{abstract}
\section{Introduction}

Let $p>3$ be a prime number and $q$ a power of p. Let $E$ be an
elliptic curve over the finite field $\mathbb{F}_q$. Thus, we can
assume that $E$ has equation $ E: y^2=x^3+ax+b$.

 The set of rational
points on $E$, denoted by $E(\mathbb{F}_q)$, has group structure. If
$n$ is an integer, we denote by $E(\mathbb{F}_q)[n]$ (or $E[n]$ if
the field is the algebraic closure $\overline{\mathbb{F}_q}$ of
$\mathbb{F}_q$) the rational points of order $n$. If $n$ is
relatively prime with $p$, $E[n] \cong \mathbb{Z}/n\mathbb{Z} \times
\mathbb{Z}/n\mathbb{Z}$.

Let $\psi_n(x)$ be the division polynomials of $E$ (see
\cite{Silverman}). As it is well known, the roots of the polynomial
$\psi_n$ are the abscissas of the $n$-torsion points, that is
$$
P=(x,y) \in E[n] \Leftrightarrow \psi_n(x)=0.
$$

Hence, the factorisation patterns of these polynomial give
information about the extension where the $n$-torsion points are
defined.

The Frobenius endomorphism,
$$
\begin{array}{cccc}
\varphi:&E(\overline{\mathbb{F}_q})&\rightarrow&E(\overline{\mathbb{F}_q})\\
&(x,y)&\rightarrow&(x^q,y^q)\\
\end{array}
$$
characterizes the rationality of a point of the elliptic curve as
follows
$$
\forall P \in E(\overline{\mathbb{F}_q}), \, P \in
E(\mathbb{F}_{q^n}) \Leftrightarrow \varphi^n(P)=P.
$$

In the paper \emph{Factorisation of division polynomials} (Proc.
Japan Academy, Ser A. 80, no. 5, pp. 79–-82), Verdure gives the
degree and the number of  factors of the division polynomial of an
elliptic curve. However, the result present there contains a
mistake. We correct it here.

\section{Patterns of $l$-th division polynomials}

Let $l$ be an odd prime different from the characteristic of
$\mathbb{F}_q$. We present here the factorisation patterns of
division polynomial only when the $l$-torsion points generate
different extension fields (the wrong result in \cite{Verdure}). If
all $l$-torsion points are defined over the same extension field,
the factorisation can be found in \cite{Verdure}.

First of all, we fix the notation. Let $f$ be a one variable
polynomial over a field $K$ of degree $n$. We say that the
factorisation pattern of $f$ is
$$
((\alpha_1,n_1),\ldots,(\alpha_d,n_d))
$$
if $f$ factorizes over $K$ as
$$
f=k\prod_{i=1}^{d}\prod_{j=1}^{n_i}P_{i,j}
$$
with $P_{i,j}$ an irreducible polynomial of degree $\alpha_i$.

The next result shows how the Frobenius endomorphism acts on $E[l]$
when the $l$-torsion points are not all defined over the same
extension of $\mathbb{F}_q$.

\begin{lem}[\cite{Verdure}]\label{actionFrob}
Let $E$ be an elliptic curve defined over $\mathbb{F}_q$. Let
$\alpha$ be the degree of the minimal extension over which an
$l$-torsion point is defined, $l$ an odd prime not equal to the
characteristic of $\mathbb{F}_q$. Assume that $E[l] \not\subset
E(\mathbb{F}_{q^\alpha})$. Then there exist $\rho \in
\mathbb{F}_l^*$ of order $\alpha$ and a basis $P,Q$ of $E[l]$ over
$\mathbb{F}_l$ in which  the n-th power of the Frobenius
endomorphism  can be expressed, for all $n$, as:

$$
\left(%
\begin{array}{cc}
  \rho^n & 0 \\
  0 & (\frac q\rho)^n \\
\end{array}%
\right) \qquad
\left(%
\begin{array}{cc}
  \rho^n & n\rho^{n-1} \\
  0 & \rho^n \\
\end{array}%
\right)
$$
if $\rho^2 \neq q$ or $\rho^2=q$ respectively.   The number $\rho$
is uniquely defined by the above properties.
\end{lem}

The previous result help us to determine the factorisation pattern
of division polynomial $\psi_l(x)$ when its factors are not all of
the same degree. The next proposition solves the mistake, in the
function $i(x,y)$, made in \cite{Verdure}.

\begin{pro}
  Let $E$ be an elliptic curve defined over $\mathbb{F}_q$. Let
  $\alpha$ be the degree of the minimal extension over which $E$ has
  a non-zero $l$-torsion point. Assume that $E[l] \not\subset
  E(\mathbb{F}_{q^\alpha})$. Let $\rho \in \mathbb{F}^*_l$ be as
  defined in Lemma \ref{actionFrob}. Let $\beta$ be the order of $q/\rho$ in
  $\mathbb{F}^*_l$. Then the pattern of the division polynomial
  $\psi_l$ is:\\

\begin{flushleft}
$
((h(\alpha),\frac{l-1}{2h(\alpha)}),(h(\beta),\frac{l-1}{2h(\beta)}),
(i(\alpha,\beta),\frac{(l-1)^2}{2i(\alpha,\beta)})) $
\end{flushleft}

\begin{flushright}
if $ q\neq \rho^2$,
\end{flushright}

\begin{flushleft}
$
((h(\alpha),\frac{l-1}{2h(\alpha)}),(h(\alpha)l,\frac{l-1}{2h(\alpha)}))
$
\end{flushleft}

\begin{flushright}
if $ q= \rho^2$,
\end{flushright}
with
$$
h(x)=\begin{cases}
x, &{\mbox{x odd,}} \\
\frac{x}{2}& {\mbox{x even,}}
\end{cases},
$$
and
$$ i(x,y)=\begin{cases}
\frac{lcm(x,y)}{2},& {\mbox{ x,y even and }} \upsilon_2(x)=\upsilon_2(y), \\
lcm(x,y),& {\mbox{otherwise.}}
\end{cases}
$$
\end{pro}

\begin{rem} Verdure gives the function $i(x,y)=~lcm(x,y)/2$ when $x$
and $y$ are both even. \end{rem}

\begin{proof}

We follow the proof given in \cite{Verdure} except for the wrong
cases.

Let $I$ be an irreducible factor $\psi_l(x)$ of degree $d$, and $P$
a point of $l$-torsion corresponding to one of its roots, then $d$
is the minimum positive integer $n$ such that $\varphi^n(P)= \pm P$.
Let $(P,Q)$ be a basis of $E[l]$ as in Lemma \ref{actionFrob}. We
distinguish the cases $q \neq \rho^2$ and $q=~\rho^2$.

\begin{itemize}
\item[i)] Suppose that $q \neq \rho^2$. If $R$ is an $l$-torsion
point which is a non-zero multiple of $P$ (or $Q$), we have that the
minimum n such that $\varphi^n(R)=\pm R$ is $n=h(\alpha)$ (or
$h(\beta)$). Notice that, $\varphi^n(R)=-R$ if and only if $\alpha$
(or $\beta$) is even, and hence $n=\alpha/2$ (or $\beta/2$).

Finally, let $R$ be any non-zero $l$-torsion point not of the
previous form, then $R=k(P+jQ)$ with $1 \leq j,k \leq l-1$. So,
$\varphi^n(R)=k(\varphi^n(P)+j\varphi^n(Q))$. The subgroup generated
by $R$ ($\langle  R\rangle$) is rational over $\mathbb{F}_{q^n}$ if
and only if $\varphi^n(R)=\pm R$. The minimum extension where
$\langle  R\rangle$ is defined is $\mathbb{F}_{q^n}$, with $n$
minimum such that $\varphi^n(R)=\pm R$.

It is easy to prove that $\varphi^n(R)=R$ if and only if
$\varphi^n(P)=P$ and $\varphi^n(Q)=Q$. Hence $lcm(\alpha,
\beta)\mid~n$ and $n=lcm(\alpha, \beta)$ is the minimum.

On the other hand, $\varphi^n(R)=-R$, if and only if
$\varphi^n(P)=-P$ and $\varphi^n(Q)=-Q$. This is only possible when
$\alpha$ and $\beta$ are both even. Moreover, $lcm(\alpha/2,
\beta/2) \mid n$ and $\alpha$ or $\beta$  not divides $lcm(\alpha/2,
\beta/2)$ (if, for example, $\alpha~\mid~lcm(\alpha/2, \beta/2)$,
then $\varphi^n(P)=P$). On the other hand, $\alpha/2$ and $\beta/2$
have the same parity, otherwise, for example, if $\alpha/2$ is even
and $\beta/2$ odd then $lcm(\alpha/2,\beta/2)= lcm(\alpha/2,\beta)$
and $\beta$ divides $lcm(\alpha/2,\beta/2)$ which is a
contradiction. If $\upsilon_2(\alpha)=\upsilon_2(\beta)$, then
$n=lcm(\alpha/2,\beta/2)$ is the minimum integer such that
$\varphi^n(P)=-P$ and $\varphi^n(Q)=-Q$. Otherwise, if both
valuations are not equal, $lcm(\alpha/2,\beta/2) $ is divisible by
$\alpha$ if $\upsilon_2(\alpha)<\upsilon_2(\beta)$ (by $\beta$ if
$\upsilon_2(\alpha)>\upsilon_2(\beta)$) which contradicts
$\varphi^n(R)=-R$.

Counting the number of points of each type, namely $l-1$, $l-1$ and
$(l-1)^2$, we have the number of factors of each type.

\item [ii)] Suppose that $q = \rho^2$. A point which is a non-zero
multiple of $P$ leads to factors of degree $\alpha$ or $\alpha/2$ as
before. If $R$ is not a multiple of $P$, then in order to have
$\varphi^n(R)=\pm R$, we have that $\rho^n=\pm 1$ and
$n\rho^{n-1}=0$. Then, depending on the parity of $\alpha$, we have
$n=lcm(\alpha,l)$ or $n=lcm(\alpha/2,l)$. Finally, since $\alpha
\mid l-1$, it is relatively prime to $l$. Therefore, these values
are $h(\alpha)l$.\end{itemize}\end{proof}

\begin{ex}
  Consider the elliptic curve $y^2=x^3+3x+6$
over $\mathbb{F}_{17}$ and take $l=5$. Then $\alpha=2$ and
$\beta=4$. According to \cite{Verdure}, the pattern of $\psi_5(x)$
should be $((1,2),(2,1),(2,4))$, but in fact it is
$((1,2),(2,1),(4,2))$.

\end{ex}
\vspace*{.5cm}

I am grateful to the referee, for this careful review and his
suggestion at the preliminar version.

\end{document}